\documentclass{article}
\usepackage{amsmath}[1996/11/01]
\usepackage{amssymb,amsthm,amsxtra}
\usepackage{eepic}

\def\[#1\]{\begin{equation}#1\end{equation}}
\makeatletter
\def\beq{%
   \relax\ifmmode
      \@badmath
   \else
      \ifvmode
         \nointerlineskip
         \makebox[.6\linewidth]%
      \fi
      $$
   \fi
}
\def\eeq{%
   \relax\ifmmode
      \ifinner
         \@badmath
      \else
         $$
      \fi
   \else
      \@badmath
   \fi
   \ignorespaces
}

\def\enddisplaymath{\eeq\global\@ignoretrue}
\makeatother

\newtheorem{thm}{Theorem}
\newtheorem{cor}[thm]{Corollary}
\newtheorem{lem}[thm]{Lemma}

\theoremstyle{remark}

\newtheorem*{rem}{Remark}
\newtheorem{rems}{Remark}[thm]

\numberwithin{thm}{section}
\numberwithin{equation}{section}

\DeclareMathOperator{\Exp}{{\bf E}}
\DeclareMathOperator{\pf}{pf}
\DeclareMathOperator{\sgn}{sgn}
\newcommand{\Z}{\mathbb Z}
\newcommand{\N}{\mathbb N}
\newcommand{\C}{\mathbb C}

\newcommand\psymmU{%
\begin{picture}(1,1)(0,0)%
\allinethickness{0.5pt}%
\path(0,0)(0,1)(1,1)(1,0)(0,0)%
\end{picture}}
\newcommand\psymmUU{%
\begin{picture}(1,1)(0,0)%
\allinethickness{0.5pt}%
\path(0,0)(0,1)(1,1)(1,0)(0,0)%
\put(0.5,0.5){\makebox(0,0){$\cdot$}}%
\end{picture}}
\newcommand\psymmO{%
\begin{picture}(1,1)(0,0)%
\allinethickness{0.5pt}%
\path(0,0)(0,1)(1,1)(1,0)(0,0)%
\path(0,0)(1,1)%
\end{picture}}
\newcommand\psymmS{%
\begin{picture}(1,1)(0,0)%
\allinethickness{0.5pt}%
\path(0,0)(0,1)(1,1)(1,0)(0,0)%
\path(1,0)(0,1)%
\end{picture}}
\newcommand\psymmu{%
\begin{picture}(1,1)(0,0)%
\allinethickness{0.5pt}%
\path(0,0)(0,1)(1,1)(1,0)(0,0)%
\path(0,0)(1,1)%
\path(0,1)(1,0)%
\end{picture}}

\newbox\tsymmUbox
\newbox\tsymmUUbox
\newbox\tsymmObox
\newbox\tsymmSbox
\newbox\tsymmubox
\setbox\tsymmUbox =\hbox{\kern0.75pt\setlength{\unitlength}{6pt}\psymmU \kern0.75pt}
\setbox\tsymmUUbox=\hbox{\kern0.75pt\setlength{\unitlength}{6pt}\psymmUU\kern0.75pt}
\setbox\tsymmObox =\hbox{\kern0.75pt\setlength{\unitlength}{6pt}\psymmO \kern0.75pt}
\setbox\tsymmSbox =\hbox{\kern0.75pt\setlength{\unitlength}{6pt}\psymmS \kern0.75pt}
\setbox\tsymmubox =\hbox{\kern0.75pt\setlength{\unitlength}{6pt}\psymmu \kern0.75pt}
\def\tsymmU{{\copy\tsymmUbox}}
\def\tsymmUU{{\copy\tsymmUUbox}}
\def\tsymmO{{\copy\tsymmObox}}
\def\tsymmS{{\copy\tsymmSbox}}
\def\tsymmu{{\copy\tsymmubox}}
\newbox\symmUbox
\newbox\symmUUbox
\newbox\symmObox
\newbox\symmSbox
\newbox\symmubox
\setbox\symmUbox =\hbox{\kern0.75pt\setlength{\unitlength}{4.5pt}\psymmU \kern0.75pt}
\setbox\symmUUbox=\hbox{\kern0.75pt\setlength{\unitlength}{4.5pt}\psymmUU\kern0.75pt}
\setbox\symmObox =\hbox{\kern0.75pt\setlength{\unitlength}{4.5pt}\psymmO \kern0.75pt}
\setbox\symmSbox =\hbox{\kern0.75pt\setlength{\unitlength}{4.5pt}\psymmS \kern0.75pt}
\setbox\symmubox =\hbox{\kern0.75pt\setlength{\unitlength}{4.5pt}\psymmu \kern0.75pt}
\def\symmU{{\copy\symmUbox}}
\def\symmUU{{\copy\symmUUbox}}
\def\symmO{{\copy\symmObox}}
\def\symmS{{\copy\symmSbox}}
\def\symmu{{\copy\symmubox}}

\begin{document}
\title{Correlation functions for symmetrized increasing subsequences}
\author{Eric M. Rains\footnote{AT\&T Research, Florham Park, New Jersey
07932, rains@research.att.com}}

\date{May 22, 2000}
\maketitle

\begin{abstract}
We show that the correlation functions associated to symmetrized increasing
subsequence problems can be expressed as pfaffians of certain antisymmetric
matrix kernels, thus generalizing the result of \cite{math.RT/9907127}
for the unsymmetrized case.
\end{abstract}

\section*{Introduction}

In \cite{math.RT/9907127}, Okounkov derived the following symmetric
function identity:  For any finite subset $S\subset \Z$,
\[
\sum_{\lambda:S\subset \{\lambda_j-j:j\in \Z^+\}}
s_\lambda(x) s_\lambda(y)
=
\sum_\lambda s_\lambda(x) s_\lambda(y)
\det(K(S)),\notag
\]
where $K(S)$ is the appropriate principal minor of an explicit infinite
matrix $K$, and $\lambda$ ranges over partitions.  The main applications of
this result are to the asymptotic analysis of generalized increasing
subsequence problems; such a problem induces a distribution on partitions
such that $\lambda$ occurs with probability $s_\lambda(x) s_\lambda(y)$,
appropriately specialized (see Section 7 of \cite{BR1}).  For instance, the
distribution of the $k$th row of $\lambda$ can be computed from this result
in terms of a certain Fredholm determinant.

In \cite{BR1}, \cite{BR2}, \cite{BR3}, we considered five classes of
generalized increasing subsequence problems, corresponding to different
choices of symmetry imposed on the problem.  As the above result only
applies to the symmetry-free class $\tsymmU$, it is natural to wonder
whether analogous results hold in the other cases.  As we shall see in
the present note, there is a matrix associated to each of the five
symmetry classes such that the corresponding correlation functions
are given as either the determinant or the pfaffian of appropriate
minors.  Each of these symmetry classes corresponds to an appropriate
Cauchy-Littlewood type identity; using the present techniques, we can
obtain analogous results for the remaining three Littlewood identities
(see Section 7).

We begin in Section 1 by giving a fairly general theorem (Theorem
\ref{thm:fpf}), inspired by the results of \cite{TracyWidom:cluster}, to
the effect that for any measure space $(X,\lambda)$ and any probability
distribution on $X^{2m}$ with density of the form
\[
\det(\phi_j(x_k)) \pf(\epsilon(x_j,x_k)),
\notag\]
the corresponding correlation function can be expressed as a pfaffian.
Since the distributions we are interested in are not of this form, we
cannot directly apply Theorem \ref{thm:fpf}.  However, in each case, we
can write the desired correlation function as a formal limit of
correlation functions to which Theorem \ref{thm:fpf} does apply.
Section 2 gives some lemmas on formal inverses of infinite matrices
which we use in sections 3 through 7 to simplify the obtained pfaffian kernels.
Finally, in section 8, we discuss the analogue for pfaffians of the notion
of Fredholm determinant, and give a Fredholm pfaffian-based derivation of
Theorem \ref{thm:fpf}.

For the (somewhat involved) definitions of the increasing subsequence
problems considered below, we refer the reader to Section 7 of
\cite{BR1}; we will also use the somewhat more general notion of
parameter set introduced in \cite{Rains:mean_identity}.

\section{Correlation functions as pfaffians}

The correlation functions we will be studying below can all be expressed as
pfaffians of certain antisymmetric matrix kernels.  Recall that a matrix
kernel on a space $X$ is a matrix-valued function on $X\times X$; for a
matrix kernel $K$, we define its transpose $K^t$ by
\[
K^t(x,y) = K(y,x)^t.
\]
Given a finite sequence $\Sigma = x_1,x_2,\dots x_k$ of elements of $X$,
the restriction $K(\Sigma)$ of $K$ to $\Sigma$ is defined to be the block
matrix with $ij$th block $K(x_i,x_j)$; note that $K^t(\Sigma)=K(\Sigma)^t$.
In particular, if $K$ is antisymmetric, then so is $K(\Sigma)$, and thus
we can compute the pfaffian $\pf(K(\Sigma))$.  When $K$ is
even-dimensional, this is invariant under reordering of $\Sigma$, and thus
depends only on the underlying set.  For a finite subset $S\subset X$,
we define $\pf(K(S))$ accordingly.  By convention, the pfaffian of
a $0\times 0$ matrix is 1, so $\pf(K(\emptyset))=1$.  Given two sequences
$\Sigma_\pm$, we define $K(\Sigma_+,\Sigma_-)$ in the obvious way, and
write $K(S_+,S_-)$ for sets $S_\pm$ whenever the meaning is clear.  Thus,
for instance, if $S_+$ and $S_-$ are disjoint, we can write
\[
\pf(K(S_+\cup S_-))
=
\pf
\pmatrix
K(S_+,S_+) & K(S_+,S_-)\\
K(S_-,S_+) & K(S_-,S_-)
\endpmatrix.
\]
We also adopt corresponding notations for determinants.

The way in which such pfaffians arise in the sequel is via the following
theorem:

\begin{thm}\label{thm:fpf}
Let $(X,\lambda)$ be a measure space, let $\phi_1,\dots \phi_{2m}$, be
functions from $X$ to $\C$, let $\epsilon$ be an antisymmetric function
from $X\times X$ to $\C$, and assume the antisymmetric matrix
\[
M_{jk}
=
\int_{x,y\in X} \phi_j(x) \epsilon(x,y) \phi_k(y) \lambda(dx) \lambda(dy)
\]
is well-defined and invertible.  For a finite subset $S=\{x_1,x_2,\dots
x_l\}\subset X$ with $l\le 2m$, we define a correlation function
\[
R(S;\phi,\epsilon)
:=
{1\over (2m-l)! \pf(M)}
\int_{x_{l+1},\dots x_{2m}\in X}
\det(\phi_j(x_k))
\pf(\epsilon(x_j,x_k))
\prod_{l+1\le j\le 2m} \lambda(dx_j);
\]
for $|S|>2m$, we set $R(S;\phi,\epsilon)=0$.  Then
$R(S;\phi,\epsilon)=\pf(K(S))$, where $K$ is the antisymmetric matrix kernel
\begin{align}
K(x,y)
=
\pmatrix
\sum_{1\le j,k\le 2m} \phi_j(x) M^{-t}_{jk} \phi_k(y)&
\sum_{1\le j,k\le 2m} \phi_j(x) M^{-t}_{jk} (\epsilon\cdot \phi_k)(y)\\
\sum_{1\le j,k\le 2m} (\epsilon\cdot \phi_j)(x) M^{-t}_{jk} \phi_k(y)&
-\epsilon(x,y)+
\sum_{1\le j,k\le 2m} (\epsilon\cdot \phi_j)(x) M^{-t}_{jk} (\epsilon\cdot \phi_k)(y)
\endpmatrix,
\end{align}
and for a function $f:X\to \C$,
\[
(\epsilon\cdot f)(x) = \int_{y\in X} \epsilon(x,y) f(y) \lambda(dy).
\]
\end{thm}

\begin{proof}
We first consider the case $|S|\ge 2m$.  In that case, if the matrix
$\Phi:=\phi_j(S)$ is singular, then the odd rows of $K(S)$ are linearly
dependent and thus $\pf(K(S))=0$.  We may thus assume $|S|=2m$ and
$\Phi$ is nonsingular.  Then we can express
$(\epsilon\cdot\phi_j)(x)$ on $S$ as a linear combination of the
functions $\phi_j(x)$.  Using this we find that
\[
\pf(K(S)) = \pf(K'(S)),
\]
where
\[
K'(x,y) = 
\pmatrix
\sum_{1\le j,k\le 2m} \phi_j(x) M^{-t}_{jk} \phi_k(y)&0\\
0&-\epsilon(x,y)
\endpmatrix.
\]
But then
\[
\pf(K'(S))
=
\pf(\Phi M^{-t} \Phi^t)\pf(\epsilon(x_j,x_k))
=
\pf(M)^{-1} \det(\phi_j(x_k)) \pf(\epsilon(x_j,x_k)),
\]
as required.

Now, suppose we know the theorem for sets of size $\ge l$, and let
$S$ be a set of size $l-1$.  Then
\[
R(S;\phi,\epsilon)
=
{1\over 2m-l+1} \int_{x_l\in X} R(S\cup \{x_l\};\phi,\epsilon) \lambda(dx_l)
=
{1\over 2m-l+1} \int_{x_l\in X} \pf(K(S\cup \{x_l\})) \lambda(dx_l)
\]
It thus suffices to show
\[
\int_{x_l\in X} \pf(K(S\cup \{x_l\})) \lambda(dx_l)
=
(2m-l+1) \pf(K(S)).
\]

Expand $\pf(K(S\cup \{x_l\}))$ along the bottom two rows and integrate,
then simplify using the following integrals:
\begin{align}
\int_{x_l\in X} K(x_l,x_l)_{21} \lambda(dx_l) &= -2m\\
\int_{x_l\in X} K(x_l,x_j)_{11} K(x_l,x_k)_{21} \lambda(dx_l)
&=
K(x_j,x_k)_{11}\\
\int_{x_l\in X} K(x_l,x_j)_{12} K(x_l,x_k)_{21} \lambda(dx_l)
&=
K(x_j,x_k)_{21}\\
\int_{x_l\in X} K(x_l,x_j)_{11} K(x_l,x_k)_{22} \lambda(dx_l)
&=
0\\
\int_{x_l\in X} K(x_l,x_j)_{12} K(x_l,x_k)_{22} \lambda(dx_l)
&=
0
\end{align}
We thus see that the $22$ terms contribute nothing.  For the $21$ terms,
$K(x_l,x_l)_{21}$ contributes $2m\pf(K(S))$ directly, while the
terms associated to $K(x_l,x_k)_{21}$ give precisely the expansion of
$\pf(K(S))$ along the first $x_k$ column, up to an overall sign change.
We thus obtain a total of $2m\pf(K(S)) - (l-1)\pf(K(S))$, as required.
\end{proof}

\begin{rems}
The above operator essentially appeared in \cite{TracyWidom:cluster}, which
considered the case $\phi_j\propto x^{j-1}$,
$\epsilon(x,y)=\frac{1}{2}\sgn(y-x)$; that reference did not obtain a direct formula for the correlation functions, however.  See Section \ref{sec:appendix}
for a derivation of the theorem along their lines.  The above proof generalizes
that used (for the same special case) in \cite{Mehta}, Chapter 6.  Note that
in \cite{Mehta}, the correlation functions are stated as ``quaternion
determinants'', essentially the restriction of the notion of pfaffian
to block matrices.
\end{rems}

\begin{rems}
When $S=\emptyset$, we find
\[
{1\over (2m)!}
\int_{x_1,\dots x_{2m}\in X}
\det(\phi_j(x_k))
\pf(\epsilon(x_j,x_k))
\prod_{1\le j\le 2m} \lambda(dx_j)
=
\pf(M),
\]
proving a result of \cite{deBruijn}.
\end{rems}

\begin{rems}
The kernel $K$ is, of course, not unique; for instance, we may use
$K'(x,y)=T(x)K(x,y)T(y)^t$ where $T$ is any function from $X$ to 
$SL_2(\C)$.
\end{rems}

\begin{cor}\label{cor:fpO}
Let $(X,\lambda)$ and $(Y,\mu)$ be measure spaces, let
$\phi_1,\dots \phi_{2m}$ be measurable functions from $X\to\C$,
let $\psi_1,\dots \psi_{2m}$ be measurable functions from
$Y\to\C$, and let $\kappa$ be a function from $X\times Y$ to $\C$.
Assume that the antisymmetric matrix
\[
M_{jk} =
\int_{x\in X, y\in Y} (\phi_j(x)\psi_k(y)-\phi_k(x)\psi_j(y)) \kappa(x,y) \lambda(dx) \mu(dy)
\]
is well-defined and invertible.  Then, for finite sets $S_0=\{x_1,x_2,\dots
x_{l_0}\}\subset X$, $S_1=\{y_1,y_2,\dots y_{l_1}\}\subset Y$, define
\begin{multline}
R(S_0,S_1;\phi,\psi,\kappa)
=
{1\over (m-l_0)!(m-l_1)!\pf(M)}
\int_{\substack{x_{l_0+1},\dots x_m\in X\\
y_{l_1+1},\dots y_m\in Y}}
\det(\phi_j(x_k)\ \psi_j(y_k))
\det(\kappa(x_j,y_k))
\\
\prod_{l_0+1\le j\le m} \lambda(dx_j)
\prod_{l_1+1\le j\le m} \mu(dy_j),
\end{multline}
we have
\[
R(S_0,S_1;\phi,\psi,\kappa)
=
\pf\pmatrix
K_{00}(S_0,S_0) & K_{01}(S_0,S_1)\\
K_{10}(S_1,S_0) & K_{11}(S_1,S_1)
\endpmatrix,
\]
where
\begin{align}
K_{00}(x,x') &=
\pmatrix
\sum_{1\le j,k\le 2m} \phi_j(x) M^{-t}_{jk} \phi_k(x')&
\sum_{1\le j,k\le 2m} \phi_j(x) M^{-t}_{jk} (\kappa\cdot \psi_k)(x')\\
\sum_{1\le j,k\le 2m} (\kappa\cdot \psi_j)(x) M^{-t}_{jk} \phi_k(x')&
\sum_{1\le j,k\le 2m} (\kappa\cdot \psi_j)(x) M^{-t}_{jk} (\kappa\cdot \psi_k)(x')
\endpmatrix
\\
K_{01}(x,y) &=
\pmatrix
\sum_{1\le j,k\le 2m} \phi_j(x) M^{-t}_{jk} (\kappa^t\cdot \phi_k)(y)&
\sum_{1\le j,k\le 2m} \phi_j(x) M^{-t}_{jk} \psi_k(y)\\
\kappa(x,y)
+\sum_{1\le j,k\le 2m} (\kappa\cdot \psi_j)(x) M^{-t}_{jk} (\kappa^t\cdot \phi_k)(y)&
\sum_{1\le j,k\le 2m} (\kappa\cdot \psi_j)(x) M^{-t}_{jk} \psi_k(y)
\endpmatrix
\\
K_{10}(y,x) &=
\pmatrix
\sum_{1\le j,k\le 2m} (\kappa^t\cdot \phi_j)(y) M^{-t}_{jk} \phi_k(x)&
-\kappa(x,y)
+\sum_{1\le j,k\le 2m} (\kappa^t\cdot \phi_j)(y) M^{-t}_{jk} (\kappa\cdot \psi_k)(x)\\
\sum_{1\le j,k\le 2m} \psi_j(y) M^{-t}_{jk} \phi_k(x)&
\sum_{1\le j,k\le 2m} \psi_j(y) M^{-t}_{jk} (\kappa\cdot \psi_k)(x)
\endpmatrix
\\
K_{11}(y,y') &=
\pmatrix
\sum_{1\le j,k\le 2m} (\kappa^t\cdot \phi_j)(y) M^{-t}_{jk} (\kappa^t\cdot \phi_k)(y')&
\sum_{1\le j,k\le 2m} (\kappa^t\cdot \phi_j)(y) M^{-t}_{jk} \psi_k(y')\\
\sum_{1\le j,k\le 2m} \psi_j(y) M^{-t}_{jk} (\kappa^t\cdot \phi_k)(y')&
\sum_{1\le j,k\le 2m} \psi_j(y) M^{-t}_{jk} \psi_k(y')
\endpmatrix
\end{align}
for $x,x'\in X$, $y,y'\in Y$.
\end{cor}

\begin{proof}
Define functions $\phi^+$ on $X\uplus Y$ by
\begin{align}
\phi^+_j(x) = \phi_j(x)\quad & \quad \phi^+_j(y) = \psi_j(y)
\end{align}
and an antisymmetric function $\epsilon$ on $(X\uplus Y)^2$ by
\begin{align}
\epsilon(x,x')=0 \quad & \quad \epsilon(x,y) = \kappa(x,y)\\
\epsilon(y,x)=-\kappa(x,y) \quad & \quad \epsilon(y,y') = 0
\end{align}
Then the function
\[
\det(\phi^+_j(z_k))
\pf(\epsilon(z_j,z_k))
\]
on $(X\uplus Y)^{2m}$ is 0 unless exactly half of the $z_k$ are in $Y$,
in which case it equals
\[
\det(\phi_j(x_k)\ \psi_j(y_k))
\det(\kappa(x_j,y_k)).
\]
Furthermore, the current matrix $M$ is the same as the matrix associated to
$\phi^+$ and $\epsilon$.  We thus find that
\[
R(S_0,S_1;\phi,\psi,\kappa)
=
R(S_0\cup S_1;\phi^+,\epsilon),
\]
so we can apply Theorem \ref{thm:fpf}; we compute
\[
(\epsilon\cdot \phi^+)(x) = (\kappa\cdot \psi)(x)\quad\quad
(\epsilon\cdot \phi^+)(y) = -(\kappa^t \cdot \phi)(y),
\]
thus obtaining the desired result, up to transformation by
\[
T(x) = \pmatrix 1&0\\0&1\endpmatrix\quad\quad
T(y) = \pmatrix 0&-1\\1&0\endpmatrix.
\]
\end{proof}

\begin{rem}
If $Y=X$, $\psi=\phi$, then we obtain a density on pairs of disjoint
$m$-subsets of $X$.  Taking the union, we obtain a density on $2m$-subsets
of $X$, which is of precisely the form considered in Theorem \ref{thm:fpf},
with $\epsilon=\kappa-\kappa^t$.  Thus the corollary may be viewed as
a refinement of the theorem, as opposed to simply a special case.
\end{rem}

\begin{cor}
Let $(X,\lambda)$ be a measure space, let $\phi_1,\dots \phi_{2m}$,
and $\psi_1,\dots \psi_{2m}$ be measurable functions from $X$ to $\C$, and
assume the antisymmetric matrix
\[
M_{jk} = \int_{x\in X} \phi_j(x)\psi_k(x)-\phi_k(x)\psi_j(x) \lambda(dx)
\]
is well-defined and invertible.  Then, defining
\[
R(S;\phi,\psi)
=
{1\over (m-l)! \pf(M)}
\int_{x_{l+1},\dots x_m\in X}
\det(\phi_j(x_k)\ \psi_j(x_k))
\prod_{l+1\le j\le m} \lambda(dx_j),
\]
we have
\[
R(S;\phi,\psi)
=
\pf(K(S)),
\]
where $K$ is the antisymmetric matrix kernel
\begin{align}
K(x,y)
&=
\pmatrix
\sum_{1\le j,k\le 2m} \phi_j(x) M^{-t}_{jk} \phi_k(y)&
\sum_{1\le j,k\le 2m} \phi_j(x) M^{-t}_{jk} \psi_k(y)\\
\sum_{1\le j,k\le 2m} \psi_j(x) M^{-t}_{jk} \phi_k(y)&
\sum_{1\le j,k\le 2m} \psi_j(x) M^{-t}_{jk} \psi_k(y)
\endpmatrix
\end{align}
\end{cor}

\begin{proof}
Apply the previous result with $(Y,\mu)=(X,\lambda)$,
$\kappa(x,y) = \delta_{xy}$, and $g=0$.
\end{proof}

In certain cases, the pfaffians simplify to determinants:

\begin{cor}\label{cor:fdu}
Let $(X,\lambda)$ and $(Y,\mu)$ be measure spaces, let
$\phi_1,\dots \phi_m$ be measurable functions from $X\to\C$,
let $\psi_1,\dots \psi_m$ be measurable functions from
$Y\to\C$, and let $\kappa$ be a function from $X\times Y\to \C$.
Assume that the matrix
\[
M_{jk} =
\int_{x\in X, y\in Y} \phi_j(x)\kappa(x,y)\psi_k(y)\lambda(dx) \mu(dy)
\]
is well-defined and invertible.  Then, defining
\begin{multline}
R_D(S_0,S_1;\phi,\psi,\kappa)
=
{1\over (m-l_0)!(m-l_1)!\det(M)}
\int_{\substack{x_{l_0+1},\dots x_m\in X\\
y_{l_1+1},\dots y_m\in Y}}
\det(\phi_j(x_k))\det(\psi_j(y_k))
\det(\kappa(x_j,y_k))
\\
\prod_{l_0+1\le j\le m} \lambda(dx_j)
\prod_{l_1+1\le j\le m} \mu(dy_j),
\end{multline}
we have
\[
R_D(S_0,S_1;\phi,\psi,\kappa)
=
\det\pmatrix
K_{00}(S_0,S_0)&K_{01}(S_0,S_1)\\K_{10}(S_1,S_0)&K_{11}(S_1,S_1)
\endpmatrix,
\]
where
\begin{align}
K_{00}(x,x') &=\sum_{1\le j,k\le m} \phi_j(x) M^{-t}_{jk} (\kappa\cdot \psi_k)(x')\\
K_{01}(x,y) &=\sum_{1\le j,k\le m} \phi_j(x) M^{-t}_{jk} \psi_k(y)\\
K_{10}(y,x) &=
-\kappa(x,y)
+\sum_{1\le j,k\le m} (\kappa^t\cdot \phi_j)(y) M^{-t}_{jk} (\kappa\cdot \psi_k)(x)\\
K_{11}(y,y') &=
\sum_{1\le j,k\le m} (\kappa^t\cdot \phi_j)(y) M^{-t}_{jk} \psi_k(y')
\end{align}
for $x,x'\in X$, $y,y'\in Y$.
\end{cor}

\begin{cor}\label{cor:fdetU}
Let $(X,\lambda)$ be a measure space, let $\phi_1,\dots
\phi_m$, and $\psi_1,\dots \psi_m$ be measurable functions from $X$ to
$\C$, and assume the matrix
\[
M_{jk} = \int_{x\in X} \phi_j(x) \psi_k(x) \lambda(dx)
\]
is well-defined and invertible.  Then, defining
\[
R_D(S;\phi,\psi)
=
{1\over (m-l)! \det(M)}
\int_{x_{l+1},\dots x_m\in X}
\det(\phi_j(x_k))
\det(\psi_j(x_k))
\prod_{l+1\le j\le m}
\lambda(dx_j),
\]
we have
\[
R_D(S;\phi,\psi)
=
\det(K(S)),
\]
where
\[
K(x,y) = \sum_{1\le j,k\le m} \phi_j(x) M^{-t}_{jk} \psi_k(y).
\]
\end{cor}

\section{Matrix inversions}

In the cases considered below, the matrices $M$ are principal minors
of certain infinite matrices; it thus becomes crucial to determine
how the inverses of the minors are related to the minors of the inverse.
The key property of the matrices is that their coefficients decay as one
gets farther away from the main diagonal.

We recall that a filtration on a ring $R$ is a sequence $R=I_0\supsetneq
I_1\supsetneq I_2\dots$ of ideals of $R$ such that $I_j I_k\subset I_{j+k}$
and $\cap_{1\le j} I_j=\{0\}$.  Equivalently, a filtration can be specified
by a valuation, that is a function $v:(R-\{0\})\to\N$ such that
\[
v(xy) \ge v(x)+v(y),\quad v(x+y)\ge \min(v(x),v(y));
\]
we simply take $v(x)=j$ whenever $I_j$ is the largest ideal in the
filtration containing $x$.  The ring $R$ is complete with respect to
the valuation $v$ if $R$ is the projective limit of the rings $R/I_j$;
equivalently, for any sequence $x_1,x_2,\dots\in R$ such that
\[
\lim_{n\to\infty} \min_{j\ne k\ge n} v(x_j-x_k) = \infty,
\]
there exists an element $x\in R$ with
\[
\lim_{n\to\infty} v(x_n-x) = \infty.
\]
The canonical example of a complete ring is a ring of formal power series,
with valuation given by the degree map.

Given an infinite matrix $M$, we let $M(m)$ denote the $m$th principal
minor of $M$.

\begin{lem}\label{lem:matU}
Let $R$ be a ring complete with respect to the valuation $v$,
and let $M$ be a matrix in $R^{\Z^+\times \Z^+}$ with decaying valuations
\[
v(M_{jk})\ge |j-k|
\]
and with unit diagonal elements.  Then $M$ is invertible,
\[
v(M^{-1}_{jk})\ge |j-k|,
\]
and for any $m\in \Z^+$,
\begin{align}
v((M(m)^{-1} - M^{-1}(m))_{jk}) &\ge 2m+2-j-k\\
v((M(m) - M^{-1}(m)^{-1})_{jk}) &\ge 2m+2-j-k.
\end{align}
In particular, for $j,k$ fixed,
\begin{align}
\lim_{m\to\infty} M(m)^{-1}_{jk} &= M^{-1}_{jk}\\
\lim_{m\to\infty} M^{-1}(m)^{-1}_{jk} &= M_{jk}.
\end{align}
\end{lem}

\begin{proof}
We first observe that for any $m$, $\det(M(m))$ is a unit in $R$; indeed,
it agrees to valuation 1 with the unit product $\prod_{1\le j\le m}
M_{jj}$.  Now, multiplication by a unit leaves the valuation unchanged, so
$v(M(m)^{-1}_{jk})=v(M(m)^{-1}_{jk}\det(M))$.  This latter element is (up
to sign) simply the determinant of the complementary minor to $(k,j)$; we
easily see that every term of this determinant has valuation at least
$m+1-j-k$.

Now, let us consider how $M(m-1)$ is related to $(M(m)^{-1})(m-1)^{-1}$.
Recall that for a block matrix
\[
M_0 = \pmatrix A&B\\C&D\endpmatrix
\]
with $D$ invertible, the upper left block of $M_0^{-1}$ is given by $(A-B
D^{-1} C)^{-1}$.  In other words, the difference between the upper left
block of $M_0$ and the inverse of the upper left block of $M_0^{-1}$ is is
$B D^{-1} C$.  Applying this to $M(m)$, we find that
\[
(M(m-1) - M(m)^{-1}(m-1)^{-1})_{jk}
={M(m)_{jm} M(m)_{mk}\over M(m)_{mm}};
\]
since $M(m)_{mm}$ is a unit, we find
\[
v((M(m-1) - M(m)^{-1}(m-1)^{-1})_{jk})
\ge
v(M(m)_{jm})+v(M(m)_{mk})
=
2m-j-k.
\]
By symmetry, we also find
\[
v((M(m)^{-1}(m-1) - M(m-1)^{-1})_{jk})
\ge
2m-j-k.
\]
By induction on $n$, we find that
\begin{align}
v((M(m) - M(n)^{-1}(m)^{-1})_{jk})&\ge 2m+2-j-k,\\
v((M(n)^{-1}(m) - M(m)^{-1})_{jk})&\ge 2m+2-j-k.
\end{align}
In particular, defining an infinite matrix $N$ by
\[
N_{jk} = \lim_{n\to\infty} M(n)^{-1}_{jk},
\]
we find $MN = NM = 1$, and the lemma follows.
\end{proof}

\begin{lem}\label{lem:matOS}
Let $R,v$ be as above, and let $M$ be an infinite antisymmetric matrix such
that
\[
v(M_{jk})\ge |j-k|-1,
\]
and $M_{(2j-1)(2j)}\in R^*$ for all $j\ge 1$.  Then $M$ is invertible and
for all $m>0$,
\begin{align}
v((M(2m) - M^{-1}(2m)^{-1})_{jk})
&\ge 4m+1-j-k.\\
v((M(2m)^{-1}-M^{-1}(2m))_{jk})
&\ge \cases 2m+2+(j+1\bmod 2)-k & k>j\\
            2m+2+(k+1\bmod 2)-j & j>k.
     \endcases
\end{align}
In particular, for $j,k$ fixed,
\begin{align}
\lim_{m\to\infty} M(m)^{-1}_{jk} &= M^{-1}_{jk}\\
\lim_{m\to\infty} M^{-1}(m)^{-1}_{jk} &= M_{jk}.
\end{align}
\end{lem}

\begin{proof}
The proof is essentially as above; the main difference is that the matrix
$D$ is now 2-dimensional, of the form
\[
D = \pmatrix 0 & u\\-u&0\endpmatrix,
\]
for some unit $u$.  Then, since $C = -B^t$, $(B D^{-1} C)_{jk}$ is
essentially just the determinant of a $2\times 2$ submatrix of $B$.
For the first equation, it is trivial to determine the valuation of this
determinant; for the second equation, we simply relate the determinant of
a $2\times 2$ minor of $M(2m)^{-1}$ to the determinant of the complementary
minor of $M(2m)$, and again the valuation is easy to determine.
\end{proof}

Similarly,

\begin{lem}\label{lem:matu}
Let $R,v$ be as above, and let $M$ be an infinite antisymmetric matrix such
that
\[
v(M_{jk})\ge |\lceil j/2\rceil - \lceil k/2\rceil|
\]
and $M_{(2j-1)(2j)}\in R^*$ for all $j\ge 1$.  Then $M$ is invertible and
for all $m>0$,
\begin{align}
v((M(2m) - M^{-1}(2m)^{-1})_{jk})
&\ge 2m+2-\lceil j/2\rceil - \lceil k/2\rceil\\
v((M(2m)^{-1}-M^{-1}(2m))_{jk})
&\ge 2m+2-\lceil j/2\rceil - \lceil k/2\rceil.
\end{align}
In particular, for $j,k$ fixed,
\begin{align}
\lim_{m\to\infty} M(m)^{-1}_{jk} &= M^{-1}_{jk}\\
\lim_{m\to\infty} M^{-1}(m)^{-1}_{jk} &= M_{jk}.
\end{align}
\end{lem}

We digress to consider a specific matrix which arises below.  For numbers
$\alpha$, $\beta$, we define $F(\alpha,\beta)$ to be the antisymmetric
matrix with
\[
F(\alpha,\beta)_{jk} = 
\cases
\alpha^{k-j-1} \beta^{(j+1)\bmod 2} \beta^{k\bmod 2}& k>j\\
-\alpha^{j-k-1} \beta^{(k+1)\bmod 2} \beta^{j\bmod 2}& j<k.
\endcases
\]
Also, if $\phi(z)$ is a Laurent series, we define the Toeplitz matrix
\[
T(\phi(z))_{jk} = [z^{k-j}] \phi(z).
\]

The following is straightforward to verify:

\begin{lem}
For any $\alpha$, $\beta\in R$ such that $v(\alpha),v(\beta)>0$,
\begin{align}
F(\alpha,\beta) &= F(-\alpha,-\beta)\\
F(\alpha,\beta)^{-1} &= -F(-\beta,\alpha),
\end{align}
and
\begin{align}
F(\alpha,1)
&= T((1-\alpha z)^{-1}) F(0,1) T((1-\alpha z)^{-1})^t\\
&= T((1-\alpha/z)^{-1}) F(0,1) T((1-\alpha/z)^{-1})^t\\
F(1,\beta)
&=
T(1+\beta z) F(1,0) T(1+\beta z)^t\\
&=
T(1+\beta/z) F(1,0) T(1+\beta/z)^t\\
F(1,0)
&=
T((1-z^2)^{-1}) F(0,1) T((1-z^2)^{-1})^t\\
&=
T((1-z^{-2})^{-1}) F(0,1) T((1-z^{-2})^{-1})^t.
\end{align}
\end{lem}

\section{The ordinary cases: $\tsymmU$ and $\tsymmUU$}

It will be instructive to rederive the result of \cite{math.RT/9907127},
since this will suggest how to deal with the symmetrized cases later.

\begin{thm}
Let $p_+$, $p_-$ be compatible parameter sets (in the sense of
\cite{Rains:mean_identity}).  Then for any finite subset
$S\subset \Z$, the probability that the set
$\{\lambda^\symmU_j(p_+,p_-)-j\}$ contains $S$ is given by
\[
\det(K^\symmU(S\mid p_+,p_-)),
\]
where
\[
K^\symmU(a,b\mid p_+,p_-)=
\sum_{1\le l}
L^\symmU(a+l\mid p_+,p_-)
L^\symmU(b+l\mid p_-,p_+)
\]
and
\[
L^\symmU(a\mid p_+,p_-) = 
[z^a] {E(z;p_+)\over E(z^{-1};p_-)},
\]
defined by contour integration over a contour containing 0 and the zeros
of $E(z^{-1};p_-)$ and excluding $\infty$ and the poles of $E(z;p_+)$.
\end{thm}

\begin{proof}
Since
\[
\Pr(\lambda^\symmU_j(p_+,p_-) = \lambda)
=
H(p_+,p_-)
s_{\lambda'}(p_+) s_{\lambda'}(p_-),
\]
we see that the theorem reduces formally to the symmetric function identity
\[
{
\sum_{\lambda:S\subset \{\lambda_i-i\}}
s_{\lambda'}(x) s_{\lambda'}(y)
\over
\sum_\lambda s_{\lambda'}(x) s_{\lambda'}(y)
}
=
\det(K^\symmU(S\mid x,y)).
\]
We first prove this formal identity, then consider the specific
specialization of interest.

If we restrict $\lambda$ so that $\ell(\lambda)\le m$, then this
only changes the left-hand-side by terms of order $O(x^m y^m)$;
it will thus suffice to derive a kernel for each $m$ such that
the formal limit $m\to\infty$ of these kernels is $K^\symmU$.

When $\ell(\lambda)\le m$, we find
\[
s_{\lambda'}(x) s_{\lambda'}(y)
=
\det(e_{\lambda_k-k+j}(x))_{j,k}\
\det(e_{\lambda_k-k+j}(y))_{j,k}.
\]
Thus we can apply Corollary \ref{cor:fdetU} above, with
\[
\phi_j(a) = e_{a+j}(x)\quad
\psi_j(a) = e_{a+j}(y).
\]
Defining $M(m)$ by
\[
M(m)_{jk} = \sum_a \phi_j(a)\psi_k(a),
\]
we find that $M(m)$ is the $m$th principal minor of the infinite
matrix
\[
M_{jk} = \sum_a e_{a+j}(x) e_{a+k}(y) = \sum_a e_{a-k}(x) e_{a-j}(y),
\]
for $1\le j,k$.  Since $j,k>0$, we can restrict the second sum to
$a>0$, and thus have
\[
M = T(E(z;y)) T(E(z;x))^t.
\]
(Recall $T(\phi(z))_{jk} = [z^{k-j}]\phi(z)$.)
We thus find
\[
M^{-1} = T(E(z;x)^{-1})^t T(E(z;y)^{-1}),
\]

With respect to the natural valuation on the ring of symmetric functions in
two variables, $M$ satisfies the hypotheses of Lemma \ref{lem:matU} above;
we thus find
\[
\lim_{m\to\infty} (M(m)^{-1} - M^{-1}(m))_{jk} = 0
\]
for any fixed $j,k$.  Since $v(\phi_j(a)) \ge a+j$, we find
\[
\lim_{m\to\infty} \sum_{1\le j,k\le m} \phi_j(a) M(m)^{-t}_{jk} \psi_k(a)
=
\sum_{1\le j,k} \phi_j(a) M^{-t}_{jk} \psi_k(a)
=
\sum_{1\le l}
(\sum_{1\le j} \phi_j(a) E(y)^{-1}_{lj})
(\sum_{1\le k} \psi_k(a) E(x)^{-1}_{lk}).
\]
We compute
\[
\sum_{1\le j} \phi_j(a) E(y)^{-1}_{lj}
=
\sum_{1\le j}
[z^{a+j}] E(z;x)
[z^{j-l}] E(z;y)^{-1}
=
\sum_j
[z^{a+j}] E(z;x)
[z^{j-l}] E(z;y)^{-1}
=
[z^{a+l}] {E(z;x)\over E(1/z;y)},
\]
thus proving the desired formal result.

For any complex number $u$ and any parameter set $p$, we define a
specialization $up$ on the ring of symmetric functions in $x$ by
\[
e_j(up) = u^j e_j(p).
\]
Now, specialize the formal identity by $e_j(x)\to e_j(up_+)$ and $e_j(y)\to
e_j(up_-)$.  For $u$ in a neighborhood of 0, both sides converge, and thus
must agree in this neighborhood.  Since both sides are analytic in a
neighborhood of the interval $[0,1]$, it follows that they must agree at
$u=1$, and the theorem is proved.
\end{proof}

\begin{rems}
Since
\[
{E(z;p_+)\over E(1/z;p_-)}
=
{H(-1/z;p_-)\over H(-z;p_+)},
\]
we find that our operator is the same as the operator of
\cite{math.RT/9907127} and \cite{math.CA/9907165} whenever the latter
operator is defined.
\end{rems}

\begin{cor}
For any finite disjoint subsets $S_+$, $S_-\subset \Z$, the probability
that the set $\{\lambda^\symmU_i(p_+,p_-)-i\}$ contains
$S_+$ and is disjoint from $S_-$ is given by
\[
\det\pmatrix
K^\symmU(S_+,S_+\mid p_+,p_-)& \sqrt{-1} K^\symmU(S_+,S_-\mid p_+,p_-)\\
\sqrt{-1} K^\symmU(S_+,S_-\mid p_+,p_-)& I-K^\symmU(S_-,S_-\mid p_+,p_-)
\endpmatrix
\]
\end{cor}

\begin{proof}
Set $T:=\{\lambda^\symmU_i(p_+,p_-)-i\}$.  Then the given determinant is
\[
\sum_{S_0\subset S_-} (-1)^{|S_0|} \Pr(S_+\cup S_0\subset T)
=
\Pr(S_+\subset T, S_-\cap T=\emptyset),
\]
as required.
\end{proof}

For the case $\tsymmUU$ of signed permutations, the analogous expectation
is a specialization of the symmetric function identity for $\tsymmU$;
we thus have:

\begin{cor}
Let $p_+$, $p_-$ be compatible parameter sets.  Then for any finite subset
$S\subset \Z$, the probability that the set
$\{\lambda^\symmUU_j(p_+,p_-)-j\}$ contains $S$ is given by
\[
\det(
K^\symmUU(S\mid p_+,p_-)
),
\]
where
\[
K^\symmUU(a,b\mid p_+,p_-)=
\sum_{1\le l}
L^\symmU((a+l)/2\mid p_+,p_-)
L^\symmU((b+l)/2\mid p_-,p_+),
\]
defining $L^\symmU(a\mid p_+,p_-):=0$ if $a\notin \Z$.
\end{cor}

\begin{proof}
After specializing, $L^\symmU(a\mid p_+,p_-)$ becomes
\[
\{[z^a] E(-z^{-2};p_-)^{-1} E(-z^2;p_+)\}
=
(-1)^{a/2} L^\symmU(a/2\mid p_+,p_-).
\]
Conjugating by $(-1)^{a/2}$ gives
\[
K^\symmUU(a,b\mid p_+,p_-)=
\sum_{1\le l}
(-1)^{a+l}
L^\symmU((a+l)/2\mid p_+,p_-)
L^\symmU((b+l)/2\mid p_-,p_+);
\]
since $L^\symmU((a+l)/2\mid p_+,p_-)=0$ unless $a+l$ is even, the result
follows.
\end{proof}

\begin{cor}
For any finite disjoint subsets $S_+$, $S_-\subset \Z$, the probability
that the set $\{\lambda^\symmUU_i(p_+,p_-)-i\}$ contains
$S_+$ and is disjoint from $S_-$ is given by
\[
\det\pmatrix
K^\symmUU(S_+,S_+\mid p_+,p_-)& \sqrt{-1} K^\symmUU(S_+,S_-\mid p_+,p_-)\\
\sqrt{-1} K^\symmUU(S_+,S_-\mid p_+,p_-)& I-K^\symmUU(S_-,S_-\mid p_+,p_-)
\endpmatrix
\]
\end{cor}

\section{The first involution case: $\tsymmO$}

Let $\delta_{a>b}$ denote the function on $\Z\times \Z$ which is 0
when $a\le b$ and $1$ when $a>b$.

\begin{thm}\label{thm:fpfO}
Let $p$ be a self-compatible parameter set, let $\alpha$ be a number with
$0\le \alpha<R(p)^{-1}$, and let $p^+$ be the parameter set obtained by
adjoining $\alpha$ to $r(p)$.  Then for any finite sets $S_0,S_1\subset \Z$,
the probability that the set $\{\lambda^\symmO_{2j-1}(p;\alpha)-2j+1\}$
contains $S_1$ and the set $\{\lambda^\symmO_{2j}(p;\alpha)-2j\}$
contains $S_0$ is given by
\[
\pf
\pmatrix
K^\symmO_{00}(S_0,S_0\mid p;\alpha)&K^\symmO_{01}(S_0,S_1\mid p;\alpha)\\
K^\symmO_{10}(S_0,S_0\mid p;\alpha)&K^\symmO_{11}(S_1,S_1\mid p;\alpha)
\endpmatrix
\]
where for $u,v\in \{0,1\}$,
\[
K^\symmO_{uv}(a,b|p;\alpha) =
\pmatrix
S^\symmO_{uv}(a,b\mid p;\alpha)&S^\symmO_{uv}(a,b+1\mid p;\alpha)\\
S^\symmO_{uv}(a+1,b\mid p;\alpha)&S^\symmO_{uv}(a+1,b+1\mid p;\alpha)
\endpmatrix
+
\cases
\delta_{b>a}\pmatrix
\alpha^{b-a} & \alpha^{b-a+1}\\
\alpha^{b-a-1}&\alpha^{b-a}\endpmatrix & uv = 01\\
-\delta_{a>b}\pmatrix
\alpha^{a-b} & \alpha^{a-b-1}\\
\alpha^{a-b+1}&\alpha^{a-b}\endpmatrix & uv = 10
\endcases
\]
with
\begin{align}
S^\symmO_{uv}(a,b\mid p;\alpha)
&=
\sum_{l>0} L^\symmO_u(a+l+1\mid p;\alpha)L^\symmO_v(b+l\mid p;\alpha)-L^\symmO_u(a+l\mid p;\alpha)L^\symmO_v(b+l+1\mid p;\alpha)\\
L^\symmO_0(a\mid p;\alpha) &= L^\symmO(a\mid p)\\
L^\symmO_1(a\mid p;\alpha) &= L^\symmO(a-1\mid p^+)\\
L^\symmO(a\mid p) &=
\delta_{\text{$a$ even}}
-
\sum_{0<j} L^\symmU(a-2j\mid p,p)
\end{align}
\end{thm}

\begin{proof}
We have
\[
\Pr(\lambda^\symmO(p;\alpha)=\lambda)
\propto
\alpha^{f(\lambda')} s_{\lambda'}(p),
\]
where $f(\lambda)$ is the number of even parts of $\lambda$, and thus
\[
\alpha^{f(\lambda')}
=
\prod_i \alpha^{\lambda_{2i-1}-\lambda_{2i-2}},
\]
so the result reduces to showing the corresponding symmetric function identity.
And again, we may take the limit $m\to\infty$ of the kernel corresponding
to the restriction $\ell(\lambda)\le 2m$.

In that case, we have
\[
\alpha^{f(\lambda')} s_{\lambda'}(x)
=
(-1)^m
\det(e_{a_k+j}(x)\ e_{b_k+j}(x))
\prod_j \alpha^{b_j-a_j-1},
\]
with
\[
a_k=\lambda_{2m-2k+2}-2m+2k-2\quad
b_k=\lambda_{2m-2k+1}-2m+2k-1.
\]
Now, if we define a kernel
\[
\kappa(a,b) = \alpha^{b-a-1} \delta_{b>a},
\]
then for nonincreasing sequences $a$ and $b$, we find
\[
\det(\kappa(a_j,b_k)) = \prod_j \alpha^{b_j-a_j-1}
\]
if $a_1<b_1\le a_2<b_2\le \dots \le a_m<b_m$; otherwise, the
determinant is 0.
We thus have
\begin{align}
\alpha^{f(\lambda')} s_{\lambda'}(x)
\propto
\det(e_{a_k+j}(x)\ e_{b_k+j}(x))
\det(\kappa(a_j,b_k))
\end{align}
for $a_1<a_2<\dots a_m$ and $b_1<b_2<\dots b_m$.  Upon symmetrizing in
$a$ and $b$, we can apply Corollary \ref{cor:fpO}, with
\[
\phi_j(a) = \psi_j(a) = e_{a+j}(x).
\]
We have
\[
(\kappa\cdot\psi_j)(a)
=
[z^{a+1+j}] (1-\alpha/z)^{-1} E(z;x),
\]
and
\[
(\kappa^t\cdot \phi_j)(a)
=
[z^{a-1+j}] (1-\alpha z)^{-1} E(z;x).
\]
Since
\begin{align}
\phi(a)+\alpha (\kappa\cdot\psi_j)(a) &= (\kappa\cdot\psi_j)(a-1)\\
\psi(a) + \alpha (\kappa^t \cdot\phi_j)(a) &= (\kappa^t \cdot \phi_j)(a+1),
\end{align}
we can simplify the matrix resulting from Corollary \ref{cor:fpO} by adding
$\alpha$ times the second row/column to the first row/column and
adding $\alpha$ times the third row/column to fourth row/column.

Now,
\[
M_{jk}
=
\sum_{a<b} (e_{a+j}(x)e_{b+k}(x)-e_{a+k}(x)e_{b+j}(x)) \alpha^{b-a-1}
=
\sum_{a<b} (e_{a-k}(x)e_{b-j}(x)-e_{a-j}(x)e_{b-k}(x)) \alpha^{b-a-1},
\]
and thus
\begin{align}
M &= T(E(z;x)) F(\alpha,1) T(E(z;x))^t\\
M^{-t} &= T(E(1/z;x)^{-1}) F(1,-\alpha) T(E(1/z;x)^{-1})^t\\
       &= T((1-\alpha/z)E(1/z;x)^{-1}(1-z^{-2})^{-1})
F(0,1)T((1-\alpha/z)E(1/z;x)^{-1}(1-z^{-2})^{-1})^t.
\end{align}
Taking $v(e_j)=j$, $v(\alpha)=1$, we see that $M$ satisfies the hypotheses
of Lemma \ref{lem:matOS} above.  Thus if $\pi,\mu$ are each either of
$\kappa\cdot \psi$, or $\kappa^t\cdot \phi$, we find
\[
\lim_{m\to\infty}
\sum_{1\le j,k\le m} \pi_j(a) M(m)^{-t}_{jk} \mu_k(b)
=
\sum_{1\le j,k} \pi_j(a) M^{-t}_{jk} \mu_k(b).
\]
It thus remains to compute
\begin{align}
\sum_{j>0}
(\kappa\cdot\psi_j)(a) T((1-\alpha/z)E(1/z;x)^{-1}(1-z^{-2})^{-1})_{jk}
&=
\sum_{j\ge 0} [z^{a+2j}] E(z;x) E(1/z;x)^{-1}\\
\sum_{j>0}
(\kappa^t\cdot\phi_j)(a) T((1-\alpha/z)E(1/z;x)^{-1}(1-z^{-2})^{-1})_{jk}
&=
\sum_{j\ge 0} [z^{a+2j}] (1-\alpha/z) E(z;x) E(1/z;x)^{-1} (1-\alpha z)^{-1}.
\end{align}
This gives the theorem, once we observe that
\[
\sum_j [z^{a+2j}] E(z;x) E(1/z;x)^{-1}
=
\delta_{\text{$a$ even}}.
\]
\end{proof}

\begin{rems}
The fact that $K_{00}$ is independent of $\alpha$ corresponds to the fact
that the joint distribution of the even rows of $\lambda^\symmO(p;\alpha)$
is independent of $\alpha$, as remarked in Section 7 of \cite{BR1}.
Similarly, the structure of $K_{11}$ corresponds to the fact that the
odd rows of $\lambda^\symmO(p;\alpha)$ are distributed as the
odd rows of $\lambda^\symmO(p^+;0)$ (which are equal to the even rows).
\end{rems}

\begin{rems}
The point of using
\[
L^\symmO(a\mid p)
=
\delta_{\text{$a$ even}} - \sum_{0<j} L^\symmU(a-2j\mid p,p)
\]
instead of
\[
L^\symmO(a\mid p)
=
\sum_{j\ge 0} L^\symmU(a+2j\mid p,p)
\]
is that the latter only converges for $p^+$ when $\alpha\le 1$ (and
converges to an incorrect value for $\alpha=1$).
\end{rems}

\begin{rems}
We observe the following relation between
$L^\symmO_0$ and $L^\symmO_1$:
\[
\alpha L^\symmO_0(a+1\mid p;\alpha)-L^\symmO_0(a\mid p;\alpha) =
\alpha L^\symmO_1(a\mid p;\alpha)-L^\symmO_1(a+1\mid p;\alpha).
\]
\end{rems}

\begin{cor}
With hypotheses as above, and $\alpha=1$, the conclusion holds with
\begin{align}
S^\symmO_{00}(a,b\mid p;1)
&=
S^\symmO_{00}(a,b\mid p)\\
S^\symmO_{01}(a,b\mid p;1)
&=
-L^\symmO(a+1\mid p)-S^\symmO_{00}(a,b\mid p)\\
S^\symmO_{10}(a,b\mid p;1)
&=
L^\symmO(b+1\mid p)-S^\symmO_{00}(a,b\mid p)\\
S^\symmO_{11}(a,b\mid p;1)
&=
L^\symmO(a+1\mid p)-L^\symmO(b+1\mid p)+S^\symmO_{00}(a,b\mid p)
\end{align}
\end{cor}

\begin{proof}
We compute
\[
{E(z;p^+)\over E(1/z;p^+)} = {-E(z;p)\over zE(1/z;p)},
\]
so
\begin{align}
L^\symmO_1(a\mid p;1) &=
\delta_{\text{$a$ odd}}
+
\sum_{0<j} [z^{a-2j}] {E(z;p)\over E(1/z;p)}\\
&=
1-L^\symmO(a\mid p).
\end{align}
\end{proof}

If we do not wish to separate the odd and even rows, we have:

\begin{cor}\label{cor:fpfO:mixed}
Let $p$ be a self-compatible parameter set, let $\alpha$ be a number with
$0\le \alpha<R(p)^{-1}$, and let $p^+$ be the parameter set obtained by
adjoining $\alpha$ to $r(p)$.  Then for any finite subset $S\subset \Z$,
the probability that $\{\lambda^\symmO_j(p;\alpha)-j\}$ contains $S$
is given by
\[
\pf(K^{\symmO'}(S\mid p;\alpha)),
\]
with
\begin{align}
K^{\symmO'}({}\mid p;\alpha)
&=
\pmatrix
S^{\symmO'}_{00}({}\mid p;\alpha)
&
S^{\symmO'}_{01}({}\mid p;\alpha)\\
S^{\symmO'}_{10}({}\mid p;\alpha)
&
S^{\symmO'}_{11}({}\mid p;\alpha)-\epsilon^\symmO({}\mid \alpha)
\endpmatrix\\
S^{\symmO'}_{uv}(a,b\mid p;\alpha) &=
\sum_{l>0} L^{\symmO'}_u(a+l+1\mid p;\alpha)L^{\symmO'}_v(b+l\mid p;\alpha)-L^{\symmO'}_u(a+l\mid p;\alpha)L^{\symmO'}_v(b+l+1\mid p;\alpha)\\
L^{\symmO'}_0(a\mid p;\alpha)
&=
(-\alpha)^{a\bmod 2} - \sum_{0<j} L^\symmU(a-2j\mid p,p^+)\\
L^{\symmO'}_1(a\mid p;\alpha)
&=
-L^\symmU(a-1\mid p^+,p)\\
\epsilon^\symmO(a,b\mid \alpha) &= \alpha^{|b-a|-1}\sgn(b-a).
\end{align}
\end{cor}

\begin{proof}
The key step is to sum over the subsets of $S$.  By the theorem, we have
\begin{align}
\sum_{S'\subset S}
\Pr(S\subset \{\lambda_j-j\})
&=
\sum_{S_0,S_1\subset S}
\pf\pmatrix
K^\symmO_{00}(S_0,S_0)&K^\symmO_{01}(S_0,S_1)\\
K^\symmO_{10}(S_1,S_0)&K^\symmO_{11}(S_0,S_1)
\endpmatrix\\
&=
\pf\left(
J+\pmatrix
K^\symmO_{00}(S,S) & K^\symmO_{01}(S,S)\\
K^\symmO_{10}(S,S) & K^\symmO_{11}(S,S)
\endpmatrix
\right),
\end{align}
where $J$ is the kernel
\[
J(a,b) = \pmatrix 0&1\\-1&0\endpmatrix.
\]

Subtract $\alpha$ times the second and third rows from the first and fourth
rows (respectively), then subtract the first row from the fourth and the
third from the second, then apply the same transformations to the columns.
This transformation is symplectic (preserves $J$), and forces the last row
of the $K$ matrix to 0.  We may thus expand along the bottom row, giving
\[
\sum_{S'\subset S}
\Pr(S\subset \{\lambda_j-j\})
=
\pf(J + K^{\symmO'}(S,S))
=
\sum_{S'\subset S} \pf(K^{\symmO'}(S')),
\]
since
\begin{align}
L^{\symmO'}_0(a\mid p;\alpha)
&=
L^\symmO_0(a\mid p;\alpha)-
\alpha L^\symmO_0(a+1\mid p;\alpha)\\
L^{\symmO'}_1(a\mid p;\alpha)
&=
L^\symmO_0(a+1\mid p;\alpha)-
L^\symmO_1(a\mid p;\alpha).
\end{align}
Thus
\[
\Pr(S\subset \{\lambda_j-j\}) = \pf(K^{\symmO'}(S'))
\]
as required.
\end{proof}

\begin{rem}
We could also have proved this directly via Theorem \ref{thm:fpf} above, with
$\phi_j(a) = e_{a+j}(x)$ and $\epsilon(a,b) = \epsilon^\symmO(a,b)$.
\end{rem}

\begin{cor}
For any finite disjoint subsets $S_+$, $S_-\subset \Z$, the
probability that $\{\lambda^\symmO_i(p;\alpha)-i\}$ contains
$S_+$ and is disjoint from $S_-$ is
\[
\pf\pmatrix K^{\symmO'}(S_+,S_+\mid p;\alpha)& \sqrt{-1} K^{\symmO'}(S_+,S_-\mid p;\alpha)\\
\sqrt{-1} K^{\symmO'}(S_-,S_+\mid p;\alpha) & J-K^{\symmO'}(S_-,S_-\mid p;\alpha)
\endpmatrix
\]
\end{cor}

\section{The second involution case: $\tsymmS$}

Similarly, for the other involution case, we have

\begin{thm}\label{thm:fpfS}
Let $p$ be a self-compatible parameter set, let $\beta$ be a number with
$0\le \beta<Q(p)^{-1}$, and let $p^+$ be the parameter set obtained
by adjoining $\beta$ to $q(p)$.  Then for any finite sets
$S_0,S_1\subset \Z$, the probability that the set
$\{\lambda^\symmS_{2j-1}(p;\beta)-2j+1\}$ contains $S_1$ and the
set $\{\lambda^\symmS_{2j}(p;\beta)-2j\}$ contains $S_0$ is given by
\[
\pf
\pmatrix
K^\symmS_{00}(S_0,S_0\mid p;\beta) &
K^\symmS_{01}(S_0,S_1\mid p;\beta) \\
K^\symmS_{10}(S_1,S_0\mid p;\beta) &
K^\symmS_{11}(S_1,S_1\mid p;\beta)
\endpmatrix,
\]
where
{
\allowdisplaybreaks
\begin{align}
K^\symmS_{00}(a,b\mid p;\beta) &=
\pmatrix
S^\symmS_{00}(a,b\mid p;\beta) &
S^\symmS_{01}(a,b\mid p;\beta) \\
S^\symmS_{10}(a,b\mid p;\beta) &
S^\symmS_{11}(a,b\mid p;\beta)
\endpmatrix
\\
K^\symmS_{01}(a,b\mid p;\beta) &=
\pmatrix
S^\symmS_{02}(a,b\mid p;\beta) &
S^\symmS_{00}(a,b\mid p;\beta) \\
S^\symmS_{12}(a,b\mid p;\beta) &
S^\symmS_{10}(a,b\mid p;\beta)
\endpmatrix
+
\delta_{b>a}
\pmatrix
0&0\\
\beta^{a\bmod 2} \beta^{(b+1)\bmod 2}&0
\endpmatrix
\\
K^\symmS_{10}(a,b\mid p;\beta) &=
\pmatrix
S^\symmS_{20}(a,b\mid p;\beta) &
S^\symmS_{21}(a,b\mid p;\beta) \\
S^\symmS_{00}(a,b\mid p;\beta) &
S^\symmS_{01}(a,b\mid p;\beta)
\endpmatrix
-
\delta_{a>b}
\pmatrix
0&\beta^{(a+1)\bmod 2} \beta^{b\bmod 2}\\
0&0
\endpmatrix
\\
K^\symmS_{11}(a,b\mid p;\beta) &=
\pmatrix
S^\symmS_{22}(a,b\mid p;\beta) &
S^\symmS_{20}(a,b\mid p;\beta) \\
S^\symmS_{02}(a,b\mid p;\beta) &
S^\symmS_{00}(a,b\mid p;\beta)
\endpmatrix
\\
S^\symmS_{uv}(a,b\mid p;\beta)
&=
\sum_{l>0} L^\symmS_u(a+l+1\mid p;\beta)L^\symmS_v(b+l\mid
p;\beta)-L^\symmS_u(a+l\mid p;\beta)L^\symmS_v(b+l+1\mid p;\beta)\\
L^\symmS_0(a\mid p;\beta)
&=
L^\symmU(a\mid p,p^+)\\
L^\symmS_1(a\mid p;\beta)
&=
\cases
\sum_{j\ge 0} L^\symmU(a+2j+1\mid p,p) & \text{$a$ even}\\
\beta \sum_{j\ge 0} L^\symmU(a+2j+2\mid p^+,p^+) & \text{$a$ odd}
\endcases\\
L^\symmS_2(a\mid p;\beta)
&=
\cases
\beta-\beta \sum_{j\ge 0} L^\symmU(a+2j\mid p,p) & \text{$a$ even}\\
1-\sum_{j\ge 0} L^\symmU(a+2j+1\mid p^+,p^+) & \text{$a$ odd}
\endcases
\end{align}
}
\end{thm}

\begin{proof}
As above, we reduce to an application of Corollary \ref{cor:fpO}, with
\[
\phi_j(a) = \psi_j(a) = e_{a+j}(x)
\]
and
\[
\kappa(a,b) = \delta_{b>a} \beta^{a\bmod 2} \beta^{(b+1)\bmod 2}.
\]
We compute
\begin{align}
(\kappa\cdot\psi_j)(a) &=
\cases
\phantom{}[z^{a+1+j}]
(1+\beta/z) E(z;x) (1-1/z^2)^{-1}
&\text{$a$ even}\\
\phantom{}[z^{a+2+j}]
\beta(1+\beta z) E(z;x) (1-1/z^2)^{-1}
&\text{$a$ odd}
\endcases
\\
(\kappa^t\cdot\phi_j)(a) &=
\cases
\phantom{}[z^{a+j}] \beta(1+\beta/z) E(z;x) z^2(1-z^2)^{-1}
&\text{$a$ even}\\
\phantom{}[z^{a+1+j}] (1+\beta z) E(z;x) z^2(1-z^2)^{-1}
&\text{$a$ odd}
\endcases
\end{align}
and
\[
M_{jk}
=
\sum_{b>a} (e_{a+j}(x) e_{b+k}(x)-e_{b+j}(x) e_{a+k}(x))
           \beta^{a\bmod 2} \beta^{(b+1)\bmod 2}.
\]
Now, when $j\bmod 2\ne k\bmod 2$, we can simply shift the variables
of summation to obtain
\[
M_{jk}
=
(T(E(z;x)) F(1,\beta) T(E(z;x))^t)_{jk}.
\]
When $j\bmod 2 = k\bmod 2$, this gives
\begin{align}
M_{jk}
&=
\sum_{b>a} 
(e_{a-k}(x) e_{b-j}(x)-e_{b-k}(x)e_{a-j}(x))
\beta^{a\bmod 2} \beta^{(b+1)\bmod 2}\\
&=
\sum_a \beta^{a\bmod 2} e_{a-j}(x)
\sum_{b>a} (e_{b-2j+k}(x)-e_{b-k}(x)) \beta^{(b+1)\bmod 2}\\
&=
\sum_a \beta^{a\bmod 2} e_{a-j}(x)
\sum_{b\le a} (e_{b-k}(x)-e_{b-2j+k}(x)) \beta^{(b+1)\bmod 2}\\
&=
\sum_{b\le a}
(e_{a-j}(x)e_{b-k}(x)-e_{a-k}(x)e_{b-j}(x))
\beta^{a\bmod 2}\beta^{(b+1)\bmod 2},
\end{align}
so we conclude that
\begin{align}
M &= T(E(z;x)) F(1,\beta) T(E(z;x))^t\\
M^{-t} &= T(E(1/z;x)^{-1}) F(-\beta,1) T(E(1/z;x)^{-1})^t\\
       &= T(E(1/z;x)^{-1}(1+\beta/z)^{-1})
          F(0,1) T(E(1/z;x)^{-1}(1+\beta/z)^{-1})^t.
\end{align}
In particular, $M^{-1}$ satisfies the hypotheses of Lemma \ref{lem:matOS},
so the kernels for finite $m$ tend to a limit.  We thus readily compute the
kernel given above.
\end{proof}

\begin{cor}\label{cor:fpfS:mixed}
Let $p$ be a self-compatible parameter set, let $\beta$ be a number with
$0\le \beta<Q(p)^{-1}$, and let $p^+$ be the parameter set obtained by
adjoining $\beta$ to $q(p)$.  Then for any finite subset $S\subset \Z$,
the probability that $\{\lambda^\symmS_j(p;\beta)-j\}$ contains $S$
is given by
\[
\pf(K^{\symmS'}(S\mid p;\beta)),
\]
with
\begin{align}
K^{\symmS'}({}\mid p;\beta)
&=
\pmatrix
S^{\symmS'}_{00}({}\mid p;\beta)
&
S^{\symmS'}_{01}({}\mid p;\beta)\\
S^{\symmS'}_{10}({}\mid p;\beta)
&
S^{\symmS'}_{11}({}\mid p;\beta)-\epsilon^\symmS({}\mid \beta)
\endpmatrix\\
S^{\symmS'}_{uv}(a,b\mid p;\beta) &=
\sum_{l>0} L^{\symmS'}_u(a+l+1\mid p;\beta)L^{\symmS'}_v(b+l\mid p;\beta)-L^{\symmS'}_u(a+l\mid p;\beta)L^{\symmS'}_v(b+l+1\mid p;\beta)\\
L^{\symmS'}_0(a\mid p;\beta)
&=
L^\symmU(a\mid p,p^+)\\
L^{\symmS'}_1(a\mid p;\beta)
&=
-\beta^{(a+1)\bmod 2}
+
\sum_{j\ge 0} L^\symmU(a+2j+1\mid p^+,p)\\
\epsilon^\symmS(a,b\mid \beta) &=
\beta^{(\max(a,b)+1)\bmod 2} \beta^{\min(a,b)\bmod 2}
\sgn(b-a).
\end{align}
\end{cor}

\section{Hyperoctahedral involutions: $\tsymmu$}

For the case $\tsymmu$ of hyperoctahedral involutions, similar arguments
can be used to derive the kernel for general $\alpha$ and $\beta$.  Since
this is rather complicated, we consider only the distribution of
$\{\lfloor\lambda_{2j-1}/2\rfloor - j\}$ and
$\{\lfloor\lambda_{2j}/2\rfloor - j\}$; or equivalently, the distribution
for $\beta=0$.

\begin{thm}\label{thm:fdu}
Let $p$ be a self-compatible parameter set, let $\alpha$ be a number with
$0\le \alpha<R(p)^{-1}$, let $\beta$ be a number with $0\le
\beta<Q(p)^{-1}$, and let $p^+$ be the parameter set obtained by adjoining
$\alpha$ to $r(p)$.  Then for any finite subsets $S_0,S_1\subset
\Z$, the probability that $\{\lfloor\lambda^\symmu_{2j-1}(p;\alpha,\beta)/2\rfloor-j\}$
contains $S_1$ and
$\{\lfloor\lambda^\symmu_{2j}(p;\alpha,\beta)/2\rfloor-j\}$ contains $S_0$ is
given by
\[
\det
\pmatrix
K^\symmu_{00}(S_0,S_0\mid p;\alpha)&K^\symmu_{01}(S_1,S_0\mid p;\alpha)\\
K^\symmu_{10}(S_0,S_1\mid p;\alpha)&K^\symmu_{11}(S_1,S_1\mid p;\alpha)
\endpmatrix,
\]
where
\begin{align}
K^\symmu_{00}(a,b\mid p;\alpha) &= \sum_{l>0} L^\symmU(a+l\mid p,p)L^\symmU(b+l\mid p,p)\\
K^\symmu_{01}(a,b\mid p;\alpha) &= \sum_{l>0} L^\symmU(a+l\mid p,p)L^\symmU(b+l\mid p,p^+)\\
K^\symmu_{10}(a,b\mid p;\alpha) &= 
\sum_{l>0} L^\symmU(a+l\mid p^+,p)L^\symmU(b+l\mid p,p)-\delta_{a\ge b} \alpha^{a-b}\\
K^\symmu_{11}(a,b\mid p;\alpha) &=
\sum_{l>0} L^\symmU(a+l\mid p^+,p)L^\symmU(b+l\mid p,p^+)
\end{align}
\end{thm}

\begin{proof}
We apply Corollary \ref{cor:fdu}, with
\[
\phi_j(a) = \psi_j(a) = e_{a+j}(x),
\]
and
\[
\kappa(a,b) = \delta_{b\ge a} \alpha^{b-a}.
\]
We find
\begin{align}
M      &= T(E(z;x)) T(E(z;x)/(1-\alpha z))^t\\
M^{-t} &= T(E(z;x)^{-1})^t T(E(z;x)^{-1}(1-\alpha z))
\end{align}
The theorem follows immediately.
\end{proof}

\begin{rem}
For general $\beta$, we instead apply Corollary \ref{cor:fpO}, with
\[
\phi_j(a) = \psi_j(a) = e_{(a+j)/2}(x)
\]
(using the convention that $e_{a/2}(x)=0$ if $a$ is odd) and
\[
\kappa(a,b) = \delta_{b\ge a} \sqrt{\alpha}^{b-a-1} \sqrt{-\beta}^{a\bmod
2}
\sqrt{-\beta}^{(b+1)\bmod 2}.
\]
We then have
\begin{align}
M &= T(E(z^2;x)) F(\sqrt{\alpha},\sqrt{-\beta}) T(E(z^2;x))^t\\
M^{-t} &= T(E(z^2;x)^{-1})^t F(-\sqrt{-\beta},\sqrt{\alpha}) T(E(z^2;x)^{-1})
\end{align}
and $M$ satisfies the hypotheses of Lemma \ref{lem:matu} above.  The
details are left to the interested reader.  (The individual terms of the
resulting operator are all fairly simple; however, since the operator
depends strongly on the parity of $a$ and $b$, there are a total of 10 such
terms to consider.)
\end{rem}

\section{Other identities}

There are three Littlewood identities that were not considered in
\cite{BR1}:
\begin{align}
\sum_{\lambda = (\alpha+1|\alpha)}
s_{\lambda'}(x)
&=
\prod_{j<k} (1+x_jx_k)\\
\sum_{\lambda = (\alpha-1|\alpha)}
s_{\lambda'}(x)
&=
\prod_j (1+x_j^2) \prod_{j<k} (1+x_jx_k)\\
\sum_{\lambda = (\alpha|\alpha)}
(-1)^{(|\lambda|-p(\lambda))/2}
s_{\lambda'}(x)
&=
\prod_j (1+x_j) \prod_{j<k} (1-x_jx_k),
\end{align}
where $(\alpha|\beta)$ is Frobenius notation, and $p((\alpha|\beta))$ is
equal to the number of parts of $\alpha$.  We also note the following
special case of the third identity:
\[
\sum_{\lambda = (\alpha|\alpha)}
\tilde{s}_{\lambda'}(x)
=
\prod_{j,k} (1+x_jx_k)
\]
For the first, second, and fourth identity, there exists an explicit
combinatorial correspondence proving the identity; in the first two
cases, this is given by \cite{Burge:correspondence}, while the
third case simply corresponds to increasing subsequences of multisets
with rotational symmetry by 90 degrees.  These correspondences extend
to the case of an arbitrary parameter set $p$ such that $p$ is compatible
with its conjugate $p'$.

As remarked in \cite{IOW}, these identities can be shown via the
Cauchy-Binet theorem.  But then Corollary \ref{cor:fdetU} implies that the
corresponding correlation functions are given in principle by appropriate
determinants.

For instance,

\begin{thm}
For any parameter set $p$ compatible with its conjugate and any
finite subset $S\subset \Z$,
\[
{
\sum_{\substack{\lambda = (\alpha-1|\alpha)\\ S\subset \{\lambda_i-i+1\}}}
s_{\lambda'}(p)
\over
\sum_{\lambda = (\alpha-1|\alpha)}
s_{\lambda'}(p)
}
=
\det(K(S)),
\]
where
\[
K(a,b)
=
(-1)^{(|b|-b)/2}
\sum_l
(-1)^{(|l|-l)/2}
L^\symmU(a+|l|\mid p,p')
L^\symmU(|b|+l\mid p',p).
\]
\end{thm}

\begin{rems}
We use $\{\lambda_i-i+1\}$ instead of $\{\lambda_i-i\}$ in order to
increase symmetry.  In particular, note that $\lambda$ is of the
appropriate form if and only if the set $\{\lambda_i-i+1\}$ contains
precisely one element of $\{j,-j\}$ for each $j$.
\end{rems}

\begin{rems}
As written, the kernel is only explicitly defined for sufficiently small
parameter sets, and must be analytically continued to the general case.
\end{rems}

\begin{proof}
For simplicity, we consider instead
\[
\sum_{\substack{\lambda = (\alpha-1|\alpha)\\ S\subset \{\lambda_i-i+1\}}}
(-1)^{|\lambda|/2}
s_{\lambda'}(p),
\]
which naturally differs only by rescaling $p$ by $\sqrt{-1}$.

We find that for $\lambda$ of the appropriate form with $\ell(\lambda)\le m$,
\[
(-1)^{|\lambda|/2}
s_{\lambda'}(p)
=
\det(\phi_j(a_k))\det(\psi_j(a_k))_{0\le j<m},
\]
where
\begin{align}
\phi_j(a) &= e_{j+a}(p)\\
\psi_j(a) &= \delta_{|a|=j}
\end{align}
and $a_k = \lambda_{m+1-i}-m+i$.  We then apply Corollary \ref{cor:fdetU}, with
\[
M_{jk} = 
\cases e_j & k=0\\
e_{j+k} + e_{j-k} & k>0
\endcases;
\]
in particular $M$ satisfies the hypotheses of Lemma \ref{lem:matU}.  We readily
verify that
\[
M^{-1}_{jk} = (-1)^j
\sum_l (-1)^l [w^{j-l} z^{k-|l|}] H(1/w;p)H(w;p)E(z;p),
\]
so
\begin{align}
\sum_{j,k\ge 0} \phi_j(a) M^{-t}_{jk} \psi_j(b)
&=
\sum_{j\ge 0} e_{j+a}(p) M^{-1}_{|b|j}\\
&=
(-1)^{|b|}
\sum_l (-1)^l [w^{|b|-l} z^{-a-|l|}] H(1/w;p)H(w;p)E(z;p)E(1/z;p)\\
&=
\sum_l [w^{|b|-l} z^{a+|l|}] {E(z;p)E(1/z;p)\over E(w;p)E(1/w;p)}.
\end{align}

Scaling $p$ by $\sqrt{-1}$ and simplifying gives the desired result.
\end{proof}

Dually,

\begin{cor}
For any parameter set $p$ compatible with its conjugate and any
finite subset $S\subset \Z$,
\[
{
\sum_{\substack{\lambda = (\alpha+1|\alpha)\\ S\subset \{\lambda_i-i+1\}}}
s_{\lambda'}(p)
\over
\sum_{\lambda = (\alpha+1|\alpha)}
s_{\lambda'}(p)
}
=
\det(I-K(S)),
\]
where
\[
K(a,b)
=
(-1)^{(|b|+b)/2}
\sum_l
(-1)^{(|l|-l)/2}
L^\symmU(-a+|l|\mid p,p')
L^\symmU(|b|+l\mid p',p).
\]
\end{cor}

For the remaining Littlewood identity, we similarly have:

\begin{thm}
For any parameter set $p$ compatible with its conjugate and any
finite subset $S\subset \Z+1/2$,
\[
{
\sum_{\substack{\lambda = (\alpha|\alpha)\\ S\subset \{\lambda_i-i+1/2\}}}
(-1)^{(|\lambda|+p(\lambda))/2}
s_{\lambda'}(p)
\over
\sum_{\lambda = (\alpha|\alpha)}
(-1)^{(|\lambda|+p(\lambda))/2}
s_{\lambda'}(p)
}
=
\det(K(S)),
\]
where
\[
K(a,b)
=
\sum_{l\in \Z+1/2} [z^{a+|l|} w^{|b|-l}] {E(z;p)E(1/z;p)\over
E(w;p)E(1/w;p)}.
\]
\end{thm}

\begin{proof}
We take
\begin{align}
\phi_j(a) &= e_{j+a+1/2}(p)\\
\psi_j(a) &= \delta_{|a|=j+1/2},
\end{align}
so
\[
M_{jk} = e_{j+k+1}(p) + e_{j-k}(p)
\]
We find
\[
M^{-1}_{jk}
=
(-1)^j
\sum_l (-1)^l [t^{j-l} u^{k+1/2-|l+1/2|}] H(1/t;p)H(t;p)E(u;p),
\]
and thus obtain the stated kernel.
\end{proof}

Specializing, we obtain (for an appropriate definition of
$\lambda^\circ(p)$, corresponding to increasing subsequences of multisets
with rotational symmetry):

\begin{cor}
Let $p$ be a parameter set compatible with its conjugate.  Then for
any finite subset $S\subset \Z+1/2$,
\[
\Pr(S\subset \{\lambda^\circ_i(p)-i+1/2\})
=
\det(K(S)),
\]
where
\[
K(a,b)
=
\sum_{l\in \Z+1/2} [z^{a+|l|} w^{|b|-l}]
{E(\sqrt{-1}z^2;p)E(\sqrt{-1}/z^2;p)\over E(\sqrt{-1}w^2;p)E(\sqrt{-1}/w^2;p)}.
\]
\end{cor}

\section{Fredholm pfaffians}\label{sec:appendix}

Let $J$ be the kernel
\[
J(a,b) = \delta_{ab} \pmatrix 0&1\\-1&0\endpmatrix
\]
Then for any other antisymmetric kernel $K$, we have
\[
\pf((J+K)(S)) = \sum_{S'\subset S} \pf(K(S')).
\]
This suggests the correct way to extend to the infinite case,
thus generalizing Fredholm determinants.  We define the
Fredholm pfaffian
\[
\pf(J+K)_X := \int_{S\subset X} \pf(K(S)) \lambda(dS),
\]
where $\lambda(dS)$ is the natural induced measure on the space of finite
subsets of $X$; by convention, $\lambda(\{\emptyset\})=1$.
In particular, when $X$ is finite and $\lambda$ is the counting measure,
we have
\[
\pf(J+K)_X = \sum_{S\subset X} \pf(K(S)) = \pf(J+K),
\]
as we would expect.  Naturally, this includes Fredholm determinants
as special cases, since
\[
\pf(J + \pmatrix \epsilon & K\\-K & 0\endpmatrix) = \det(I+K),
\]
for any scalar kernel $K$ and any antisymmetric scalar kernel $\epsilon$.

We note the following properties of Fredholm pfaffians:

\begin{lem}
For any antisymmetric matrix kernel $K$,
\[
\pf(J+K)_X^2
=
\det(I+J^{-1} K)_X.
\]
For any ordinary matrix kernel $K_0$,
\[
\pf((I+K_0) (J+K) (I+K_0^t))_X
=
\det(I+K_0)_X \pf(J+K)_X.
\]
If $A$ is a matrix operator from $X$ to $Y$, $M_X$
is an invertible antisymmetric matrix operator on $X$, and $M_Y$ is an
invertible antisymmetric matrix operator on $Y$, then
\[
\pf(M_Y)_Y
\pf(M_Y^{-t} + A M_X A^t)_Y
=
\pf(M_X)_X
\pf(M_X^{-t} + A^t M_Y A)_X.
\]
\end{lem}

\begin{rem}
The last equation generalizes the Fredholm determinant identity
\[
\det(M_1)\det(M_1^{-1}+A M_2 B)
=
\det(M_2)\det(M_2^{-1}+B M_1 A).
\]
\end{rem}

The significance of Fredholm pfaffians for our purposes is related to
the following result:

\begin{thm}\label{thm:correls}
Let $(X,\lambda)$ be a measure space, and let $\mu$ be a measure
on the set of countable subsets of $X$.  Suppose
\[
\int_{T\subset X} \chi_T(dS) \mu(dT) = \pf(K(S)) \lambda(dS),
\]
where $\chi_T$ is the atomic measure concentrated on the finite subsets
of $T$.  Then for functions $f:X\to \C$,
\[
\int_{T\subset X}
\prod_{x\in T} (1+f(x))
\mu(dT)
=
\pf(J + \sqrt{f} K \sqrt{f})_{X,\lambda}
\]
whenever both sides are defined.
\end{thm}

\begin{proof}
On the one hand, we have
\[
\int_{T\subset X}
\prod_{x\in T} (1+f(x))
\mu(dT)
=
\int_{T\subset X}
\sum_{S\subset T} \prod_{x\in S} f(x)
\mu(dT)
=
\int_{T\subset X}
\int_{S\subset X} \prod_{x\in S} f(x) \chi_T(dS) \mu(dT);
\]
on the other hand, we have
\[
\pf(J + \sqrt{f} K \sqrt{f})_{X,\lambda}
=
\int_{S\subset X} \pf(\sqrt{f} K\sqrt{f})(S) \lambda(dS)
=
\int_{S\subset X} \prod_{x\in S} f(x) \pf(K(S)) \lambda(dS).
\]
The theorem follows.
\end{proof}

\begin{rems}
Note that
\[
\pf(J + \sqrt{f} K \sqrt{f})_{X,\lambda}
=
\pf(J + K)_{X,f\lambda},
\]
where $(f\lambda)(dx) = f(x) \lambda(dx)$; thus the square root is best
thought of as merely notational.
\end{rems}

Note in particular that if $X=\Z$, $\lambda$ is the counting measure, and
$\mu$ is a probability measure, then $\Exp(\chi_T(\{S\}))$ is precisely
equal to $\Pr(S\subset T)$, thus explaining the connection with our earlier
results.

In particular, Theorem \ref{thm:fpf} is related to a Fredholm pfaffian
result:

\begin{thm}
Let $(X,\lambda)$ be a measure space, let $f$, $\phi_1,\dots \phi_{2m}$, be
functions from $X$ to $\C$, let $\epsilon$ be an antisymmetric function
from $X\times X$ to $\C$, and assume the antisymmetric matrix
\[
M_{jk}
=
\int_{x,y\in X} \phi_j(x) \epsilon(x,y) \phi_k(y) \lambda(dx) \lambda(dy)
\]
is well-defined and invertible.  Then
\begin{align}
F(f;\phi,\epsilon)
&:=
{1\over (2m)! \pf(M)}
\int_{x_1,\dots x_{2m}\in X}
\det(\phi_j(x_k))
\pf(\epsilon(x_j,x_k))
\prod_{1\le j\le 2m} (1+f(x_j)) \lambda(dx_j)\\
&\phantom{:}=
\pf(J+\sqrt{f} K\sqrt{f})_{X,\lambda},
\end{align}
in the sense that if either side is defined, then both are defined,
and take the same value.
\end{thm}

\begin{proof}
This of course follows immediately from Theorem \ref{thm:fpf}, but the
following independent proof (based on the arguments of
\cite{TracyWidom:cluster}) gives useful insight into how the kernel
$K$ can be derived.  (The above proof, of course, has the advantage of
using only finite methods.)

From Section 4 of \cite{deBruijn}, we have
\[
\int_{x_1,x_2,\dots x_{2m}}
\det(\phi_j(x_k))
\pf(\epsilon(x_j,x_k))
\prod_j \mu(dx_j)
=
(2m)!
\pf(
\int_{x,y} \phi_j(x) \epsilon(x,y) \phi_k(y)
\mu(dx) \mu(dy)
)
\]
for any measure $\mu$.  Thus, taking $\mu = (1+f)\lambda$, we find
\begin{align}
F(f;\phi,\epsilon)
&=
\pf(M)^{-1}
\pf(
\int_{x,y\in X}
\phi_j(x) \epsilon(x,y) \phi_k(y) (1+f(x))(1+f(y)) \lambda(dx) \lambda(dy)
)\\
&=
\pf(M)^{-1}
\pf(M + A M_X A^t)\\
&=
\pf(M_X)_X
\pf(M_X^{-t} + A^t M^{-t} A)_X
\end{align}
where
\[
A =
\pmatrix
\sqrt{f}\phi_j &
\sqrt{f}\epsilon\cdot\phi_j
\endpmatrix\quad\quad
M_X = 
\pmatrix
\sqrt{f}(x)\epsilon(x,y)\sqrt{f}(y) & I\\
-I & 0
\endpmatrix
\]
We thus find $\pf(M_X)_X = 1$ and
\[
M_X^{-t} = 
\pmatrix
0 & I\\
-I & -\sqrt{f}(x)\epsilon(x,y)\sqrt{f}(y)
\endpmatrix.
\]
Thus
\[
F(f;\phi,\epsilon)
=
\pf(J + \sqrt{f} K \sqrt{f})_X
\]
as required.
\end{proof}

Let $\lambda$ be a random partition.  We say that the distribution of
$\lambda$ is represented by the antisymmetric kernel $K(a,b)$ on $\Z$ if
\[
\Pr(S\subset \{\lambda_i-i\}) = \pf(K(S)).
\]
(Thus, for instance, $\lambda^\symmU(p_+,p_-)$ is represented by
\[
\pmatrix 0&K^\symmU({}\mid p_+,p_-)\\-(K^\symmU)^t({}\mid
p_+,p_-)&0\endpmatrix,
\]
and similarly for the other partition distributions considered above.)
We observe that for any set $N$, the Fredholm pfaffian
\[
\pf(J-\sqrt{t} K \sqrt{t})_N
\]
encodes the distribution of $|\{\lambda_i-i\}\cap N\}|$, and thus
as $n$ varies,
\[
\pf(J-\sqrt{t} K \sqrt{t})_{\{n,n+1,\dots\}}
\]
encodes the marginal distribution of $\lambda_i$ for each $i$.
With this in mind, we give the following Fredholm pfaffian identity:

\begin{thm}\label{thm:disccont}
Let $K$ be an antisymmetric matrix kernel that represents a probability
distribution on the set of partitions.  Then for any decomposition
$\Z = N_+\uplus N_-$ such that $N_{+-}:=N_+\cap \Z^-$ and $N_{-+}:=N_-\cap
\N$ are both finite,
\begin{align}
\pf(J-t^{1/4}(K-\chi_{N_-} J \chi_{N_-})t^{1/4})_\Z
&=
(1+\sqrt{t})^{|N_{-+}|-|N_{+-}|} \pf(J-\sqrt{t} K\sqrt{t})_{N_+},\\
&=
(1-\sqrt{t})^{|N_{+-}|-|N_{-+}|} \pf(J-\sqrt{t} (J-K)\sqrt{t})_{N_-}.
\end{align}
where $\chi_{N_-}$ is the projection onto $N_-$.
\end{thm}

\begin{proof}
Let $\lambda$ be the random partition associated to $K$, and set
$T:=\{\lambda_j-j\}$, $T_+ = T\cap N_+$, $T_- = N_--T$.
By the definition of the Fredholm pfaffian,
{
\allowdisplaybreaks
\begin{align}
\pf(J - t^{1/4}(K - \chi_{N_-} J \chi_{N_-})t^{1/4})
&=
\sum_{S\subset \Z}
t^{|S|/2}
\pf((\chi_{N_-} J \chi_{N_-} - K)(S))\\
&=
\sum_{S_\pm\subset N_\pm}
t^{(|S_+|+|S_-|)/2}
\pf\pmatrix -K(S_+,S_+) &-K(S_+,S_-)\\
            -K(S_-,S_+) &(J-K)(S_-,S_-)\endpmatrix\\
&=
\sum_{S_\pm\subset N_\pm}
t^{(|S_+|+|S_-|)/2}
(-1)^{|S_+|}
\pf\pmatrix K(S_+,S_+) &\sqrt{-1} K(S_+,S_-)\\
            \sqrt{-1} K(S_-,S_+) &(J-K)(S_-,S_-)\endpmatrix\\
&=
\sum_{S_\pm\subset N_\pm}
t^{(|S_+|+|S_-|)/2}
(-1)^{|S_+|}
\Pr(S_+\subset T, S_-\cap T=\emptyset)\\
&=
\sum_{S_\pm\subset N_\pm}
t^{(|S_+|+|S_-|)/2}
(-1)^{|S_+|}
\Pr(S_\pm\subset T_\pm)\\
&=
\sum_{R_\pm\subset N_\pm}
\sum_{S_\pm\subset R_\pm}
t^{(|S_+|+|S_-|)/2}
(-1)^{|S_+|}
\Pr(T_\pm=R_\pm)\\
&=
\sum_{R_\pm\subset N_\pm}
(1+\sqrt{t})^{|R_-|}
(1-\sqrt{t})^{|R_+|}
\Pr(T_\pm=R_\pm)\\
&=
\sum_{R_\pm\subset N_\pm}
(1-t)^{|R_+|}
(1+\sqrt{t})^{|R_-|-|R_+|}
\Pr(T_\pm=R_\pm).
\end{align}
}

Now, we have the following lemma:

\begin{lem}
Let $\Z=N_+\cap N_-$ be a decomposition as above.  Then for
any partition $\lambda$ with associated set $T$,
\[
|N_+\cap T|-|N_- - T| = |N_{+-}|-|N_{-+}|.
\]
\end{lem}

\begin{proof}
Recall that for any partition,
\[
|T\cap \N| = |\Z^- - T|.
\]
Setting $N_{++}=N_+\cap \N$, $N_{--}=N_-\cap \Z^-$, we have
\[
|T\cap \N| = |N_{++}\cap T| + |N_{-+}\cap T|
           = |N_{++}\cap T| + |N_{-+}| - |N_{-+}-T|
\]
and
\[
|\Z^- - T| = |N_{+-} - T| + |N_{--} - T|
           = |N_{+-}| - |N_{+-}\cap T| + |N_{--}-T|.
\]
Subtracting these two quantities, we conclude that
\[
|N_+\cap T| + |N_{-+}| - |N_--T| - |N_{+-}|
=
0.
\]
\end{proof}

We may thus replace $(1+\sqrt{t})^{|R_-|-|R_+|}$ in the above sum with
$(1+\sqrt{t})^{|N_{-+}|-|N_{+-}|}$.  We thus have
\begin{align}
\pf(J + t^{1/4}(K - \chi_{N_-} J \chi_{N_-})t^{1/4})
&=
(1+\sqrt{t})^{|N_{-+}|-|N_{+-}|}
\sum_{R_\pm\subset N_\pm}
(1-t)^{|R_+|}
\Pr(T_\pm=R_\pm)\\
&=
(1+\sqrt{t})^{|N_{-+}|-|N_{+-}|}
\sum_{R_+\subset N_+}
(1-t)^{|R_+|}
\Pr(T_+=R_+)\\
&=
(1+\sqrt{t})^{|N_{-+}|-|N_{+-}|}
\sum_{S_+\subset N_+}
(-t)^{|S_+|}
\Pr(S_+\subset T)\\
&=
\pf(J-\sqrt{t} K\sqrt{t})_{N_+}.
\end{align}

Similarly,
\begin{align}
\pf(J + t^{1/4}(K - \chi_{N_-} J \chi_{N_-})t^{1/4})
&=
(1-\sqrt{t})^{|N_{+-}|-|N_{-+}|}
\sum_{R_-\subset N_-}
(1-t)^{|R_-|}
\Pr(T_-=R_-)\\
&=
(1-\sqrt{t})^{|N_{+-}|-|N_{-+}|}
\sum_{S_-\subset N_-}
(-t)^{|S_-|}
\Pr(S_-\cap T = \emptyset)\\
&=
(1-\sqrt{t})^{|N_{+-}|-|N_{-+}|}
\pf(J-\sqrt{t} (J-K)\sqrt{t})_{N_-}.
\end{align}
\end{proof}

\begin{rems}
The point of the theorem is that while
\[
\pf(J + t^{1/4}(K - \chi_{N_-} J \chi_{N_-})t^{1/4})_\Z
\]
is rather more complicated as a pfaffian on $\Z$, its image under the
Fourier transform (which as an orthogonal transformation preserves
Fredholm pfaffians) is much more likely than
\[
\pf(J + t^{1/2} K t^{1/2})_{N_+}
\]
to have a simple kernel on the unit circle.  Indeed, for the first pfaffian
to have a simple kernel, all that is necessary is for $K$ and $\chi_{N_-}$
to have simple kernels; for the second pfaffian, their composition must
also be simple.
\end{rems}

\begin{rems}
Note that in particular,
\[
\pf(J-\sqrt{t} (J-K)\sqrt{t})_{N_-}
=
(1-t)^{|N_{-+}|-|N_{+-}|} \pf(J-\sqrt{t} K\sqrt{t})_{N_+}.
\]
\end{rems}

\begin{cor}
Let $K$ be a scalar kernel such that
\[
\pmatrix 0&K\\-K^t&0\endpmatrix
\]
represents a probability distribution on the set of partitions.  Then for
any decomposition $\Z=N_+\uplus N_-$ such that $N_{+-}:=N_+\cap \Z^-$ and
$N_{-+}:=N_-\cap \N$ are both finite,
\begin{align}
\det(I - t^{1/2}(K - \chi_{N_-}))_\Z
&=
(1+\sqrt{t})^{|N_{-+}|-|N_{+-}|} \det(I-t K)_{N_+},\\
&=
(1-\sqrt{t})^{|N_{+-}|-|N_{-+}|} \det(I-t (I-K))_{N_-}.
\end{align}
\end{cor}

For instance, taking $K=K^\symmU(\mid|p_+,p_-)$ and conjugating by the
Fourier transform, we find
\[
\det(1-\lambda K)_{[n,\infty)}
=
(1+\sqrt{\lambda})^{-n}
\det(I - \lambda^{1/2} K')_C,
\]
where
\begin{align}
K'(z,w) &= {z^{-n} w^n - \phi(z)\phi(w)^{-1}\over 2\pi i(z-w)},\\
\phi(z) &= {E(z;p_+)\over E(1/z;p_-)},
\end{align}
and with $C$ an appropriately chosen contour containing 0.  This
generalizes the results of \cite{BDJ2} (which essentially showed that when
$p_+=p_-=t{:}/$, the identity holds to second order at $\lambda=1$).  For a
direct, analytic proof of this identity, see \cite{BDR}.

We close by remarking that \cite{math.CA/9907165} used the identity of
\cite{math.RT/9907127} to express a large class of Toeplitz determinants as
discrete Fredholm determinants, or equivalently, to so express a large
class of integrals over the unitary group.  Similarly, Corollaries
\ref{cor:fpfO:mixed} and \ref{cor:fpfS:mixed} can be used to express
appropriate integrals over the orthogonal and symplectic groups as discrete
Fredholm pfaffians:
\begin{align}
\int_{U\in O(l)} \det(E(U;p))
&=
Z^\symmO(p;0)^{-1}
\pf(J-K^{\symmO'}({}\mid p;0))_{[l,\infty)}\\
\int_{U\in Sp(2l)} \det(E(U;p))
&=
Z^\symmS(p;0)^{-1}
\pf(J-K^{\symmS'}({}\mid p;0))_{[2l,\infty)}
\end{align}
(actually statements about {\it formal} integrals); here
\begin{align}
Z^\symmO(p;0) &:= \pf(J-K^{\symmO'}({}\mid p;0))_{[0,\infty)}\\
Z^\symmS(p;0) &:= \pf(J-K^{\symmS'}({}\mid p;0))_{[0,\infty)}.
\end{align}
We can also use Theorem \ref{thm:disccont} to rewrite these as continuous
Fredholm pfaffians; details are left to the reader.

\end{document}